\documentclass[english]{amsart}
\usepackage{babel}

\makeatletter

\providecommand{\LyX}{L\kern-.1667em\lower.25em\hbox{Y}\kern-.125emX\@}

 \theoremstyle{plain}    
 \newtheorem{thm}{Theorem}[section]
 \numberwithin{equation}{section} 
 \numberwithin{figure}{section} 
 \theoremstyle{remark}
 \newtheorem{rem}[thm]{Remark}
 \theoremstyle{remark}    
 \newtheorem*{acknowledgement*}{Acknowledgement} 

\makeatother
\begin{document}

\title{ Quasi-Lie bialgebroids and Twisted Poisson Manifolds.}

\author{Dmitry Roytenberg}

\begin{abstract}
We develop a theory of quasi-Lie bialgebroids using a homological
approach. This notion is a generalization of quasi-Lie bialgebras,
as well as twisted Poisson structures with a 3-form background which
have recently appeared in the context of string theory, and were studied
by \v{S}evera and Weinstein using a different method.
\end{abstract}
\maketitle

\section{\label{sec:Intro}Introduction.}

The purpose of this note is twofold: to develop a theory of quasi-Lie
bialgebroids and to use it to treat twisted Poisson manifolds with
a closed 3-form background. The latter appeared recently in works
of Park \cite{Park}, Cornalba-Schiappa \cite{CoSchia} and 
Klim\v{c}ik-Strobl
\cite{KlimStrobl} on string theory, and were studied by \v{S}evera
and Weinstein \cite{SevWe} using Courant algebroids and Dirac structures.
In the present work we present an alternative approach along the lines
of \cite{KS3} and \cite{Roy1}; it can be viewed as a companion to
\cite{SevWe}. 

The notion of a \emph{Lie bialgebroid} was introduced by Mackenzie
and Xu in \cite{MacXu} where it appeared as linearization of Poisson
groupoids. It consists of a pair of Lie algebroid structures on dual
vector bundles \( A \) and \( A^{*} \) satisfying a compatibility
condition. In \cite{Roy1} we developed an approach to the theory
of Lie bialgebroids in terms of a pair of mutually commuting Hamiltonians
on the symplectic supermanifold \( {\mathcal{E}}=T^{*}\Pi A=T^{*}\Pi A^{*} \).
The algebra of functions on \( {\mathcal{E}} \) has a natural double
grading, and the properties of Lie bialgebroids follow easily from
an interplay of elementary symplectic geometry and homological algebra.
In particular, the sum of the two Hamiltonians is self-commuting,
hence it defines a homological (hamiltonian) vector field \( D \)
on \( {\mathcal{E}} \); the pair \( ({\mathcal{E}},D) \) is called
the \emph{homological double} of the Lie bialgebroid \( (A,A^{*}) \).
The Courant algebroid structure on \( A\oplus A^{*} \) originally
proposed by Liu \emph{et al.} \cite{LWX1} as a {}``Drinfeld double''
of \( (A,A^{*}) \) is recovered from the homological double by Kosmann-Schwarzbach's
derived bracket construction \cite{KS5}. This generalizes the homological
approach to Lie bialgebras pioneered by Lecomte and Roger \cite{LecRog},
in which case it reduces to pure algebra.

In the same work \cite{Roy1} we remarked that a skew-symmetric trilinear
form on \( A \) or \( A^{*} \) (viewed as a function on \( {\mathcal{E}} \))
can be added to the sum of the above two Hamiltonians; requiring the
total sum to self-commute leads to the notion of a \emph{quasi-Lie
bialgebroid}, generalizing Kosmann-Schwarzbach's quasi-Lie bialgebras
treated in \cite{KS3} using the same homological approach. We observed
that a quasi-Lie bialgebroid also has a homological double which produces
a Courant algebroid structure on \( A\oplus A^{*} \). Aside from
these observations, the theory of quasi-Lie bialgebroids lay dormant
for a while due to a lack of interesting examples. The situation changed
several months ago when twisted Poisson manifolds were brought to
our attention by Alan Weinstein. These are manifolds equipped with
a bivector field \( \pi  \) and a closed 3-form \( \phi  \) (the
background field) satisfying the equation \begin{equation}
\label{eqn:twistedPoisson}
\frac{1}{2}[\pi ,\pi ]=\wedge ^{3}\tilde{\pi }\phi .
\end{equation}
 Here, and elsewhere in this note, given a bilinear form \( B \)
on a vector space \( V \), we shall denote by \( \tilde{B} \) the
corresponding map from \( V \) to \( V^{*} \)defined by \( \tilde{B}(\alpha )(\beta )=B(\alpha ,\beta ) \)
for \( \alpha ,\beta \in V \). An example of a twisted Poisson manifold
is provided by a Lie group \( G \) equipped with a bi-invariant metric
\( <\cdot ,\cdot > \). In this case \( \phi  \) is the canonical
bi-invariant Cartan 3-form \( \phi (u,v,w)=\frac{1}{2}<[u,v],w> \),
while \( \pi  \) is determined by \( \tilde{\pi }(g)=2(Ad_{g}-1)/(Ad_{g}+1) \),
where we identify both the tangent and cotangent space at \( g \)
with the Lie algebra \( {\mathfrak {g}} \) via left translation and
\( <\cdot ,\cdot > \). \( \pi  \) is defined on the open dense subset
of \( G \) where \( Ad_{g}+1 \) is invertible (see \cite{SevWe}
for details). 

The realization that twisted Poisson manifolds provide examples of
quasi-Lie bialgebroids gave impetus for developing the general theory
which resulted in the present note. The note is organized as follows.
In Section \ref{sec:bialg} we briefly recall the notion of Lie bialgebroid
emphasizing the approach of \cite{Roy1}. In Section \ref{sec:quasi}
we define a quasi-Lie bialgebroid structure on \( (A,A^{*}) \) as
a triple of Hamiltonians \( \mu  \), \( \gamma  \) and \( \phi  \)
on \( {\mathcal{E}} \) of respective degrees \( (1,2) \), \( (2,1) \)
and \( (0,3) \), such that their sum \( \Theta  \) Poisson-commutes
with itself. We unravel the resulting algebra by showing that such
structures correspond to \emph{differential quasi-Gerstanhaber algebra}
structures on \( \Gamma (\bigwedge ^{\cdot }A) \) (a term coined
by Huebschmann \cite{Huebsch-Quasi} to denote a special case of homotopy
Gerstenhaber algebras), or to \emph{quasi-differential Gerstenhaber
algebra} structures on \( \Gamma (\bigwedge ^{\cdot }A^{*}) \). We
define the homological double as \( ({\mathcal{E}},D=\{\Theta ,\cdot \}) \)
and remark on the existence of a spectral sequence converging to the
cohomology of \( D \).

In Section \ref{sec:twisting} we study the phenomenon of \emph{twisting}
(introduced in \cite{Dr2} for quasi-Hopf algebras and in \cite{KS3}
for quasi-Lie bialgebras). Here again the symplectic geometry of \( {\mathcal{E}} \)
yields a clear understanding of the phenomenon and allows for a quick
and easy derivation of formulas. Specifically, the twisting by an
\( \omega \in \Gamma (\wedge ^{2}A^{*}) \) or by \( \pi \in \Gamma (\wedge ^{2}A) \)
(thought of as functions on \( {\mathcal{E}} \)) is the canonical
transformation given by the flow of the Hamiltonian vector field \( X_{\omega }=\{\omega ,\cdot \} \)
(resp. \( X_{\pi }=\{\pi ,\cdot \} \)). It is immediately seen that
the twisting by \( \omega  \) transforms a quasi-Lie bialgebroid
structure on \( (A,A^{*}) \) into a new one (with isomorphic double),
whereas the twisting by \( \pi  \) yields a quasi-Lie bialgebroid
provided \( \pi  \) satisfies the twisted (non-homogeneous) Maurer-Cartan
equation. In case of a Lie bialgebroid, the Maurer-Cartan equation
is homogeneous; it was obtained in \cite{LWX1} by a different method.

In Section \ref{sec:symmetries} we study the action of the natural
group of symmetries: canonical transformations of \( {\mathcal{E}} \)
preserving the grading, and particularly the subgroup \( {\mathcal{G}} \)
that acts trivially on \( M\subset {\mathcal{E}} \), referred to
as the gauge group. We obtain a factorization of \( {\mathcal{G}} \)
similar to that used in the theory of Poisson-Lie groups. Also in
this context we briefly mention deformation theory governed by the
quasi-Gerstenhaber algebra associated to a quasi-Lie bialgebroid;
the appropriate structure equation is the twisted Maurer-Cartan equation.

Lastly, in Section \ref{sec:examples} we apply the theory to two
special cases: arbitrary bivector fields and twisted Poisson manifolds
with a 3-form background. We take \( A \) to be the tangent bundle
of a manifold \( M \). In the former case we start with the standard
Lie bialgebroid structure on \( (TM,T^{*}M) \) and twist by a bivector
field \( \pi  \). The result is a quasi-Lie bialgebroid structure
on \( (T^{*}M,TM) \) which is a Lie bialgebroid if and only if \( [\pi ,\pi ]=0 \).
In the latter case, we start by adding a closed 3-form \( \phi  \)
to the standard structure on \( (TM,T^{*}M) \), and then twist by
\( \pi  \) (this is the way Park \cite{Park} arrived at twisted
Poisson manifolds in string-theoretic context). The result is again
a quasi-Lie bialgebroid provided the twisted Maurer-Cartan equation
holds, which in this case reduces to (\ref{eqn:twistedPoisson}).
This extends the observation of \v{S}evera and Weinstein that a twisted
Poisson manifold defines a Lie algebroid structure on \( T^{*}M \).
The formulas obtained in \cite{SevWe} are recovered as a special
case of the general formulas from Section \ref{sec:twisting}. We
also interpret the local action of the group of 2-forms introduced
in \cite{SevWe} in terms of the factorization of the gauge group
\( {\mathcal{G}} \).

One issue that we have avoided in this note is quantization. Park
\cite{Park} argues on physical grounds that deformation quantization
of \( ({\mathcal{E}},\Theta ) \) is related to Deligne's conjecture
(now a theorem of Kontsevich) on the formality of 2-algebras; BV quantization
of the corresponding topological field theory would produce explicit
formulas (Park also considers higher \( p \)-algebras). These issues
are beyond the scope of the present note and will be approached elsewhere.

\begin{acknowledgement*}
We would like to thank Anton Alekseev, Yvette Kosmann-Schwarzbach,
Jae-Suk Park, Pavol \v{S}evera, Thomas Strobl, Alan Weinstein
and Ping Xu for useful discussions and advice. We also thank the Erwin
Schr{\"o}dinger Institute, Penn State University, MSRI and UC Berkeley
for hospitality during the time when this work was carried out.
\end{acknowledgement*}

\section{\label{sec:bialg}Lie bialgebroids and their doubles.}

Given a vector bundle \( A\rightarrow M \), \( M \) a manifold,
we denote by \( \Pi A \) the supermanifold whose algebra of functions
is \( \Gamma (\bigwedge ^{\cdot }A^{*}) \). Of course, \( \Pi A \)
is also a vector bundle over \( M \) (in fact, sometimes the notation
\( A[1] \) is used instead of \( \Pi A \) to emphasize the nonnegative
integer grading coming from the bundle structure, rather than just
the parity shift). 

A \emph{Lie algebroid} structure on \( A \) is a homological vector
field on \( \Pi A \) of degree \( +1 \), i.e. a derivation \( d_{A} \)
of \( \Gamma (\bigwedge ^{\cdot }A^{*}) \) increasing degrees by
\( 1 \) and satisfying \( [d_{A},d_{A}]=2d^{2}_{A}=0 \). Alternatively
(and equivalently), one can define a Lie algebroid structure on \( A \)
as an odd Poisson bracket of degree \( -1 \) on \( C^{\infty }(\Pi A^{*})=\Gamma (\bigwedge ^{\cdot }A) \),
usually referred to as a \emph{Schouten} or \emph{Gerstenhaber} bracket.
We denote this bracket by \( [\cdot ,\cdot ]_{A} \). The usual definition
in terms of an anchor and bracket on sections of \( A \) is recovered
easily from either of these. 

We say that \( (A,A^{*}) \) is a \emph{Lie bialgebroid} if both \( A \)
and \( A^{*} \) are Lie algebroids and in addition, the vector field
\( d_{A} \) preserves (is a derivation of) \( [\cdot ,\cdot ]_{A^{*}} \).
If \( A \) is a Lie algebroid, \( (A,A^{*}) \) becomes a Lie bialgebroid
in a trivial way if we endow \( A^{*} \) with the zero algebroid
structure. A good nontrivial example of a Lie bialgebroid is \( (T^{*}M,TM) \),
where \( M \) is a Poisson manifold. Here \( [\cdot ,\cdot ]_{A^{*}}=[\cdot ,\cdot ] \),
the usual Schouten bracket of multivector fields, while \( d_{A}=d_{\pi }=[\pi ,\cdot ] \),
where \( \pi \in \Gamma (\bigwedge ^{2}TM) \) is the Poisson tensor.
The fact that \( d^{2}_{\pi }=0 \) is equivalent to \( [\pi ,\pi ]=0 \);
the compatibility condition is obvious. The anchor on \( T^{*}M \)
is \( \tilde{\pi } \), while the bracket of 1-forms is the \emph{Koszul
bracket:} \begin{equation}
\label{eqn:koszulbracket}
[\alpha ,\beta ]_{\pi }={\mathcal{L}}_{\tilde{\pi }\alpha }\beta -{\mathcal{L}}_{\tilde{\pi }\beta }\alpha -d(\pi (\alpha ,\beta )).
\end{equation}

The following construction gives a completely symmetric picture of
Lie bialgebroids and also enables one to define the {}``double''
of a Lie bialgebroid, generalizing the Drinfeld double of a Lie bialgebra.
Given a vector bundle \( A\rightarrow M \), let \( {\mathcal{E}}=T^{*}\Pi A \).
It is an even symplectic supermanifold. Local coordinates \( \{x^{i}\} \)
on \( M \) and a local basis \( \{e_{a}\} \) of sections of \( A \)
give rise to Darboux coordinates \( \{x^{i},\xi ^{a},p_{i},\theta _{a}\} \)
on \( {\mathcal{E}} \). The symplectic form is \( \Omega =dx^{i}dp_{i}+d\xi ^{a}d\theta _{a} \).
There is a canonical symplectomorphism, the \emph{Legendre transformation}
\( L:T^{*}\Pi A\rightarrow T^{*}\Pi A^{*} \) \cite{Roy1}; in the
above coordinates, it simply exchanges \( \xi ^{a} \) with \( \theta _{a} \).
Thus, \( {\mathcal{E}} \) fibres as the cotangent bundle over both
\( \Pi A \) and \( \Pi A^{*} \). It fits into the following diagram:

\begin{equation}
\label{eqn:doublebundle}
\begin{array}{ccc}
{\mathcal{E}} & \stackrel{\bar{h}}{\longrightarrow } & \Pi A^{*}\\
\left. \bar{v}\right\downarrow  &  & \left. v\right\downarrow \\
\Pi A & \stackrel{{h}}{\longrightarrow } & M
\end{array}
\end{equation}

This diagram is an example of a \emph{double vector bundle} \cite{Prad-1,Mac-Double1,KonUrb}
i.e. each arrow is a vector bundle, and both the horizontal and the
vertical pairs of arrows are vector bundle morphisms. The double vector
bundle structure gives rise to a double (nonnegative integer) grading
of functions on the supermanifold \( {\mathcal{E}} \). Each grading
separately is not compatible with parity, but their sum, the total
weight, is. Whenever we speak of functions on \( {\mathcal{E}} \),
we shall mean polynomials with respect to the total grading; we denote
this algebra by \( {\mathcal{C}^{\cdot }} \), or \( {\mathcal{C}^{\cdot ,\cdot }} \)
when we want to emphasize the double grading. Notice that the canonical
symplectic form on \( {\mathcal{E}} \) is of bi-degree \( (1,1) \),
and thus of total weight \( 2 \); the canonical Poisson bracket has
bi-degree \( (-1,-1) \) and total weight \( -2 \). The \emph{core}
of \( {\mathcal{E}} \) is \( \ker \bar{h}\bigcap \ker \bar{v} \),
where \( \bar{h} \) (resp. \( \bar{v} \)) is viewed as a vector
bundle morphism over \( h \) (resp.\( v \)). It is a vector bundle
over \( M \), although the whole of \( {\mathcal{E}} \) is not.
In the present case the core is easily seen to be \( T^{*}M \) (more
precisely, \( T^{*}[2]M \) if the weight is taken into account);
it is also the support of the supermanifold \( {\mathcal{E}} \). 

We can also consider \( E=A\oplus A^{*} \) and the corresponding
supermanifold \( \Pi E \). It is an even Poisson manifold; the fibres
of \( \Pi E \) are the symplectic leaves, and the basic functions
are Casimir. Since \( \Pi E \) is the fibre product of \( \Pi A \)
and \( \Pi A^{*} \) over \( M \), by the universal property there
is a canonical fibration \( p:{\mathcal{E}}\rightarrow \Pi E \).
This fibration is easily seen to be a Poisson map; in fact, \( {\mathcal{E}} \)
is a minimal symplectic realization of \( \Pi E \) in the sense that
the extra dimension is equal to the co-rank of the Poisson tensor
on \( \Pi E \). The corresponding doubly graded Poisson subalgebra
of functions is \( {\bar{\mathcal{C}}}^{\cdot ,\cdot }=\Gamma (\bigwedge ^{\cdot }A\otimes \bigwedge ^{\cdot }A^{*}) \),
with \( {\bar{\mathcal{C}}}^{\cdot }=\Gamma (\bigwedge ^{\cdot }E^{*}) \).
The Poisson bracket on \( {\bar{\mathcal{C}}}^{\cdot } \) is the
{}``big bracket'' \cite{LecRog,KS3} applied pointwise over \( M \):
it is the unique extension of the canonical pairing of \( A \) and
\( A^{*} \) to a bi-derivation of the exterior algebra of sections
of \( E \) (essentially, just contraction of tensors). Notice that
\( {\bar{\mathcal{C}}}^{0,\cdot }={\mathcal{C}}^{0.\cdot } \), and
\( {\bar{\mathcal{C}}}^{\cdot ,0}={\mathcal{C}}^{\cdot ,0} \). If
\( M \) is a point, \( {\mathcal{E}} \) and \( \Pi E \) are the
same.

Notice finally that both \( \Pi A \) and \( \Pi A^{*} \) sit inside
\( {\mathcal{E}} \) as the zero sections of, respectively, \( \bar{v} \)
and \( \bar{h} \); they are, of course, Lagrangian submanifolds of
\( {\mathcal{E}} \). We denote these distinguished Lagrangian submanifolds
by \( L \) and \( L^{*} \), respectively. More generally, by elementary
symplectic geometry, \( \bar{v} \)-projectable graded Lagrangian
submanifolds are of the form \( d\omega (\Pi A) \), where \( d\omega  \)
is the exterior derivative of a quadratic function \( \omega \in \Gamma (\wedge ^{2}A^{*}) \),
viewed as a section of \( \bar{v} \). Similarly, \( \bar{h} \)-projectable
graded Lagrangian submanifolds are of the form \( d\pi (\Pi A^{*}) \)
for some \( \pi \in \Gamma (\wedge ^{2}A) \). Given \( \omega  \)
(resp. \( \pi  \)), we denote the corresponding Lagrangian submanifold
\( L_{\omega } \) (resp. \( L^{*}_{\pi } \)). The projection \( p:{\mathcal{E}}\rightarrow \Pi E \)
maps \( L_{\omega } \) (resp. \( L^{*}_{\pi } \)) onto the graph
of \( \omega  \) (resp. \( \pi  \)) in \( A\oplus A^{*} \). In
general, graded Lagrangian submanifolds that contain \( M \) correspond
under \( p \) to maximally isotropic subbundles of \( E \) that
are not necessarily graphs.

Now, a Lie algebroid structure on \( A \) is equivalent to a function
\( \mu  \) on \( {\mathcal{E}} \) of bi-degree \( (1,2) \) satisfying 

\begin{equation}
\label{eqn:algebroid}
\{\mu ,\mu \}=0
\end{equation}
 The homological vector field on \( \Pi A \) and the Schouten bracket
on \( \Pi A^{*} \) are recovered as follows. For a function \( \alpha  \)
on \( \Pi A \), 

\[
\bar{v}^{*}(d_{\mu }\alpha )=\{\mu ,\bar{v}^{*}\alpha \},\]
 and for two functions \( X,Y \) on \( \Pi A^{*} \), 

\[
\bar{h}^{*}([X,Y]_{\mu })=\{\{\bar{h}^{*}X,\mu \},\bar{h}^{*}Y\}\},\]
 the so-called \emph{derived bracket}; as \( \{\bar{h}^{*}X,\bar{h}^{*}Y\}=0 \)
\( \forall X \) and \( Y \), the derived bracket is a graded Lie
algebra bracket (for details see \cite{KS5}, where the notion was
first introduced and studied). In what follows we shall suppress the
pullback notation and write \( \alpha  \) instead of \( \bar{v}^{*}\alpha  \)
and so on, when the meaning is clear.

Likewise, a Lie algebroid structure on \( A^{*} \) is equivalent
to a function \( \gamma  \) on \( {\mathcal{E}} \) of bi-degree
\( (2,1) \) satisfying 

\begin{equation}
\label{eqn:coalgebroid}
\{\gamma ,\gamma \}=0
\end{equation}

It gives a differential \( d_{\gamma } \) on \( \Pi A^{*} \) and
a Schouten bracket \( [\cdot ,\cdot ]_{\gamma } \) on \( \Pi A \)
by analogous formulas. Moreover, it was proved in \cite{Roy1}%
\footnote{This was independently proved by A. Vaintrob (unpublished)
} that \( (A,A^{*}) \) is a Lie bialgebroid if and only if, in addition, 

\begin{equation}
\label{eqn:cocycle}
\{\mu ,\gamma \}=0
\end{equation}

It follows that the notion of Lie bialgebroid is self-dual, i.e. \( (A,A^{*}) \)
is a Lie bialgebroid if and only if \( (A^{*},A) \) is; furthermore,
\( \Theta =\mu +\gamma  \) (of total weight \( 3 \)) satisfies 

\begin{equation}
\label{eqn:structure}
\{\Theta ,\Theta \}=0
\end{equation}
if and only if (\ref{eqn:algebroid}), (\ref{eqn:coalgebroid}) and
(\ref{eqn:cocycle}) hold. Therefore, in this case its Hamiltonian
vector field \( D=\{\Theta ,\cdot \} \) is homological of degree
\( +1 \) on \( {\mathcal{E}} \). The pair \( ({\mathcal{E}},D) \)
is called the \emph{homological double} of the Lie bialgebroid \( (A,A^{*}) \). 

\begin{rem}
In view of (\ref{eqn:algebroid}), the Hamiltonian vector field \( \{\mu ,\cdot \} \)
is homological of bi-degree \( (0,+1) \). Thus, for each fixed \( k\geq 0 \),
\( ({\mathcal{C}^{k,\cdot }},\{\mu ,\cdot \}) \) is a differential
complex. For \( k=0 \) it is just the standard complex of the Lie
algebroid \( A \) with coefficients in the trivial representation;
for higher \( k \) it should be regarded as the standard complex
of \( A \) with coefficients in the \( k \)-th exterior power of
the adjoint representation, to which it reduces when the base \( M \)
is a point. In particular, (\ref{eqn:cocycle}) implies that \( \gamma  \)
is a 1-cocycle on \( A \) with coefficients in the exterior square
of the adjoint, which is exactly the way in which Lie bialgebras were
originally defined. Similarly, \( ({\mathcal{C}^{\cdot ,l}},\{\gamma ,\cdot \}) \)
is a complex for each \( l\geq 0 \), and due to (\ref{eqn:cocycle})
\( ({\mathcal{C}^{\cdot ,\cdot }},\{\mu ,\cdot \},\{\gamma ,\cdot \}) \)
is a double complex whose total complex is \( ({\mathcal{C}^{\cdot }},D=\{\Theta ,\cdot \}) \).
For Lie bialgebras, this double complex and its spectral sequence
were studied in \cite{LecRog}. For Lie bialgebroids (e.g. Poisson
manifolds), investigating the spectral sequence could yield interesting
results. It will be carried out in a separate paper.
\begin{rem}
We can observe, by counting degrees, that the space of functions on
\( {\mathcal{E}} \) of degree \( \leq 1 \) is closed under both
the Poisson bracket \( \{\cdot ,\cdot \} \) and the derived bracket
\( \{\{\cdot ,\Theta \},\cdot \} \). These functions correspond under
the projection \( p \) to the basic and linear functions on \( \Pi E \),
i.e. smooth functions on \( M \) and sections of \( E=A\oplus A^{*} \).
The Poisson bracket restricts to the canonical inner product on \( A\oplus A^{*} \),
and if \( M \) is a point, the derived bracket of linear functions
is precisely Drinfeld's double Lie bracket. In general, the derived
bracket is not skew-symmetric; the resulting structure on \( E \)
coincides (after skew-symmetrization) with the \emph{Courant algebroid}
introduced by Liu, Weinstein and Xu \cite{LWX1} as a candidate for the 
double
of \( (A,A^{*}) \); the Courant algebroid axioms follow easily from
(\ref{eqn:structure}) (see \cite{Roy1}) for details). The space
of all functions on \( \Pi E \) is not closed under the derived bracket,
but on functions on \( {\mathcal{E}} \) we get a \emph{Loday-Gerstenhaber
algebra} structure. 
\end{rem}
\end{rem}

\section{\label{sec:quasi}Quasi-Lie bialgebroids and quasi-Gerstenhaber algebras.}

The most general cubic function \( \Theta  \) on \( {\mathcal{E}} \)
is of the form \( \Theta =\mu +\gamma +\phi +\psi  \); it contains
terms \( \phi  \) of degree \( (0,3) \) and \( \psi  \) of degree
\( (3,0) \), in addition to terms \( \mu  \) and \( \gamma  \)
considered above (then \( \phi \in \Gamma (\bigwedge ^{3}A^{*}) \),
while \( \psi \in \Gamma (\bigwedge ^{3}A) \)). If \( \Theta  \)
obeys the structure equation (\ref{eqn:structure}), we get a so-called
\emph{proto Lie bialgebroid} structure on \( (A,A^{*}) \) (\cite{Roy1},
generalizing \cite{KS3}). Splitting the quartic \( \{\Theta ,\Theta \} \)
into components according to the double grading, we get the following
set of equations: \begin{equation}
\label{eqn:protoLie}
\begin{array}{rcl}
\frac{1}{2}\{\mu ,\mu \}+\{\gamma ,\phi \} & = & 0\\
\frac{1}{2}\{\gamma ,\gamma \}+\{\mu ,\psi \} & = & 0\\
\{\mu ,\gamma \}+\{\phi ,\psi \} & = & 0\\
\{\mu ,\phi \} & = & 0\\
\{\gamma ,\psi \} & = & 0
\end{array}
\end{equation}
This structure is rather complicated, but it still has a homological
double \( ({\mathcal{E}},D) \) which produces on \( E=A\oplus A^{*} \)
a Courant algebroid structure (see \cite{Roy1}). We shall mainly
restrict our attention to the case when at least one of the additional
terms (say, \( \psi  \)) vanishes. Then the equations above reduce
to:

\begin{equation}
\label{eqn:quasiLie}
\begin{array}{rcl}
\frac{1}{2}\{\mu ,\mu \}+\{\gamma ,\phi \} & = & 0\\
\{\gamma ,\gamma \} & = & 0\\
\{\mu ,\gamma \} & = & 0\\
\{\mu ,\phi \} & = & 0
\end{array}
\end{equation}
 By definition, the pair of dual vector bundles \( (A,A^{*}) \) together
with \( \mu  \), \( \gamma  \) and \( \phi  \) of degree \( (1,2) \),
\( (2,1) \) and \( (0,3) \), respectively, satisfying (\ref{eqn:quasiLie}),
is a \emph{quasi-Lie bialgebroid}; we call \( ({\mathcal{E}},D=\{\Theta ,\cdot \}) \)
its \emph{homological double}. 

Notice that quasi-Lie bialgebroids can be characterized as those self-commuting
cubic Hamiltonians on \( {\mathcal{E}} \) that vanish on the distinguished
graded Lagrangian submanifold \( L^{*} \), the image of the zero
section \( \Pi A^{*}\rightarrow T^{*}\Pi A^{*} \). Thus, a quasi-Lie
bialgebroid can be equivalently described by a \emph{Manin pair} \( (({\mathcal{E}},\Theta ),L^{*}) \).
Lie bialgebroids correspond to those \( \Theta  \) that also vanish
on the other distinguished Lagrangian submanifold \( L \) (the image
of \( \Pi A\rightarrow T^{*}\Pi A \)), and so can be described by
\emph{Manin triples} \( (({\mathcal{E}},\Theta ),L,L^{*}) \). 

Let us now unravel the resulting algebraic structure. First of all,
\( \gamma  \) defines a Lie algebroid structure on \( A^{*} \) ,
i.e. a derivation \( d=d_{\gamma } \) of \( \Gamma (\bigwedge ^{\cdot }A) \)
of degree \( +1 \) and square zero (due to \( \{\gamma ,\gamma \}=0 \)).
On the other hand, \( \mu  \) induces a bracket \( [\cdot ,\cdot ]=[\cdot ,\cdot ]_{\mu } \)
on \( \Gamma (\bigwedge ^{\cdot }A) \) of degree \( -1 \), completely
determined by an anchor map \( a:A\to TM \) and a bracket on sections
of \( A \). The compatibility condition \( \{\mu ,\gamma \}=0 \)
between \( \mu  \) and \( \gamma  \) means that \( d \) is a derivation
of \( [\cdot ,\cdot ] \). The equation \( \frac{1}{2}\{\mu ,\mu \}+\{\gamma ,\phi \}=0 \)
means that the usual Lie algebroid relations are satisfied only up
to certain {}``defects'' depeding on \( \gamma  \) and \( \phi  \).
Specifically, for all \( X,Y\in \Gamma (A) \), \begin{equation}
\label{eqn:L2L2anchors}
a([X,Y])=[a(X),a(Y)]+a_{*}\phi (X,Y)
\end{equation}
 where \( a_{*} \) is the anchor of the Lie algebroid \( A^{*} \)
and \( \phi (X,Y)=i_{X\wedge Y}\phi \in \Gamma (A^{*}) \), and for
all \( X,Y,Z\in \Gamma (A) \),  \begin{equation}
\label{eqn:L2L2brackets}
\begin{array}{l}
[[X,Y],Z]+[[Y,Z],X]+[[Z,X],Y]=\\
=d(\phi (X,Y,Z))+\phi (dX,Y,Z)-\phi (X,dY,Z)+\phi (X,Y,dZ)
\end{array}
\end{equation}
where \( \phi  \) is viewed as a bundle map \( \bigwedge ^{4}A\rightarrow A \).
Finally, \( \{\mu ,\phi \}=0 \) is a coherence condition between
\( [\cdot ,\cdot ] \) and \( \phi  \): for all \( X,Y,Z,W\in \Gamma (A) \), 

\begin{equation}
\label{eqn:L2L3}
\begin{array}{c}
([\phi (X,Y,Z),W]-[\phi (X,Y,W),Z]+[\phi (X,Z,W),Y]-[\phi (Y,Z,W),X])-\\
-(\phi ([X,Y],Z,W)-\phi ([X,Z],Y,W)+\phi ([X,W],Y,Z)+\\
+\phi ([Y,Z],X,W)-\phi ([Y,W],X,Z)+\phi ([Z,W],X,Y))=0
\end{array}
\end{equation}
 Let us define, for each positive integer \( k \), a \( k \)-linear
map \( l_{k} \) on \( \Gamma (\bigwedge ^{\cdot }A)[1] \) with values
in itself, of degree \( 2-k \), as follows. Let \( l_{1}=d \), \( l_{2}=[\cdot ,\cdot ] \),
\( l_{3}=\phi  \) (more precisely, we let \( l_{3} \) vanish if
one of the arguments is a function, let \( l_{3}(X,Y,Z)=\phi (X,Y,Z) \)
for \( X,Y,Z\in \Gamma (A) \), and then extend it as a derivation
in each argument), and \( l_{k}=0 \) for \( k>3 \). Then the above
relations imply that \( (\Gamma (\bigwedge ^{\cdot }A)[1],\{l_{k}\}) \)
is a \emph{strongly homotopy Lie} (or \( L_{\infty } \)) algebra
\cite{SHLA}. Moreover, since each \( l_{k} \) is a \( k \)-derivation
of the exterior multiplication, what we get is in fact a \emph{strongly
homotopy Gerstenhaber} (or \( G_{\infty } \)) algebra. Such \( G_{\infty } \)-algebras
are very simple compared to the most general possible case: the exterior
multiplication remains undeformed and each \( l_{k} \) is a strict
derivation. They were called \emph{differential quasi-Gerstenhaber
algebras} by Huebschmann \cite{Huebsch-Quasi}. It is easy to see
that, conversely, any \( G_{\infty } \)-algebra on \( \Gamma (\bigwedge ^{\cdot }A) \)
satisfying these conditions must come from a quasi-Lie bialgebroid.
To summarize: quasi-Lie bialgebroid structures on \( (A,A^{*}) \)
are in 1-1 correspondence with quasi-Gerstenhaber algebra structures
on \( \Gamma (\bigwedge ^{\cdot }A) \).

Dually, a quasi-Lie bialgebroid structure on \( (A,A^{*}) \) is also
equivalent to a \emph{quasi-differential Gerstenhaber algebra} structure
on \( \Gamma (\bigwedge ^{\cdot }A^{*}) \): the Gerstenhaber bracket
\( [\cdot ,\cdot ]_{\gamma } \) is given by \( \gamma  \) (it satisfies
graded Jacobi strictly), while \( \mu  \) gives a derivation (quasi-differential)
\( d_{\mu } \) of \( [\cdot ,\cdot ]_{\gamma } \) of degree \( +1 \),
which does not square to \( 0 \) but rather satisfies \[
d^{2}_{\mu }+[\phi ,\cdot ]_{\gamma }=0\]
 Futhermore, \( d_{\mu }\phi =0 \).

\begin{rem}
The equations (\ref{eqn:quasiLie}) imply that \( d_{\mu }=\{\mu ,\cdot \} \)
induces a differential on the cohomology of \( ({\mathcal{C}^{\cdot ,\cdot }},d_{\gamma }) \)
and in fact, there still exists a spectral sequence converging to
the cohomology of \( ({\mathcal{C}^{\cdot }},D=\{\Theta ,\cdot \}) \).
\end{rem}
Of course if we replace a \( \phi  \) of bi-degree \( (0,3) \) by
a \( \psi  \) of bi-degree \( (3,0) \), the roles of \( \mu  \)
and \( \gamma  \) are reversed: we get a differential quasi-Gerstenhaber
algebra on \( \Gamma (\bigwedge ^{\cdot }A^{*}) \) and a dual quasi-differential
Gerstenhaber algebra on \( \Gamma (\bigwedge ^{\cdot }A) \), i.e.
a quasi-Lie bialgebroid structure on \( (A^{*},A) \). When the base
\( M \) is a point, we recover various quasi-bialgebras studied by
Kosmann-Schwarzbach \cite{KS3}.

\section{\label{sec:twisting}Fibre translations and twisting.}

Let \( \omega \in \Gamma (\bigwedge ^{2}A^{*}) \). When pulled back
to \( {\mathcal{E}} \), \( \omega  \) has bi-degree \( (0,2) \).
Its Hamiltonian vector field \( X_{\omega }=\{\omega ,\cdot \} \)
is thus of bi-degree \( (-1,1) \). Hence, the action of \( X_{\omega } \)
on \( {\mathcal{E}} \) preserves the total weight but not the double
grading. The flow of \( X_{\omega } \) is the fibre translation along
\( -d\omega  \) with respect to the fibration \( \bar{v}:{\mathcal{E}}=T^{*}\Pi A\rightarrow \Pi A \).
Let \( F_{\omega } \) be the time \( 1 \) map of the flow. The corresponding
pullback of functions can be expressed as \[
F^{*}_{\omega }=\exp X_{\omega }=1+X_{\omega }+\frac{1}{2}X^{2}_{\omega }+\frac{1}{6}X^{3}_{\omega }+\cdots \]
 In coordinates, \( \omega =\frac{1}{2}\omega _{ab}(x)\xi ^{a}\xi ^{b} \),
and \( F_{\omega } \) is given by: \begin{equation}
\label{eqn:translbyomega}
\begin{array}{ccl}
\tilde{x}^{i} & = & x^{i}\\
\tilde{\xi }^{a} & = & \xi ^{a}\\
\tilde{p}_{i} & = & p_{i}-\frac{1}{2}\frac{\partial \omega _{ab}}{\partial x^{i}}\xi ^{a}\xi ^{b}\\
\tilde{\theta }_{a} & = & \theta _{a}-\omega _{ab}\xi ^{b}
\end{array}
\end{equation}

We are going to apply this fibre translation to a cubic hamiltonian
\( \Theta =\mu +\gamma +\phi +\psi  \). Let \( \Theta _{\omega }=F^{*}_{\omega }\Theta =(\exp X_{\omega })\Theta  \);
then \( \Theta _{\omega } \) is again cubic, and it satisfies (\ref{eqn:structure})
if and only if \( \Theta  \) does (since the flow acts by canonical
transformations). We say that \( \Theta _{\omega } \) is the \emph{twisting}
of \( \Theta  \) by \( \omega  \). Notice that the above exponential
series gets truncated when applied to a function of finite degree.
In fact, \( \Theta _{\omega }=\mu _{\omega }+\gamma _{\omega }+\phi _{\omega }+\psi _{\omega } \)
where \begin{equation}
\label{eqn:twistbyomega}
\begin{array}{ccl}
\mu _{\omega } & = & \mu +X_{\omega }\gamma +\frac{1}{2}X^{2}_{\omega }\psi =\mu +h_{[\omega ,\cdot ]_{\gamma }}+\wedge ^{2}\tilde{\omega }\psi \\
\gamma _{\omega } & = & \gamma +X_{\omega }\psi =\gamma +\tilde{\omega }\psi \\
\phi _{\omega } & = & \phi +X_{\omega }\mu +\frac{1}{2}X^{2}_{\omega }\gamma +\frac{1}{6}X^{3}_{\omega }\psi =\phi -d_{\mu }\omega -\frac{1}{2}[\omega ,\omega ]_{\gamma }+\wedge ^{3}\tilde{\omega }\psi \\
\psi _{\omega } & = & \psi 
\end{array}
\end{equation}
 Here \( \tilde{\omega }:A\rightarrow A^{*} \) is used to lower the
specified number of indices on \( \psi  \), whereas \( h_{v} \)
denotes the linear hamiltonian on \( {\mathcal{E}}=T^{*}\Pi A \)
corresponding to a vector field \( v \) on \( \Pi A \). These formulas
follow easily from the definitions.

Notice that for a quasi-Lie bialgebroid structure on \( (A,A^{*}) \)
(i.e. when \( \psi =0 \)) the twisting by \( \omega  \) automatically
produces a new quasi-Lie bialgebroid structure. On the other hand,
for a quasi-Lie bialgebroid structure on \( (A^{*},A) \) (i.e. when
\( \phi =0 \)), the twisting yields a quasi-Lie bialgebroid provided
\( \omega  \) satisfies a certain integrability condition, the \emph{twisted
Maurer-Cartan equation}: \begin{equation}
\label{eqn:twistedMC1}
d_{\mu }\omega +\frac{1}{2}[\omega ,\omega ]_{\gamma }=\wedge ^{3}\tilde{\omega }\psi 
\end{equation}

We can also twist by elements \( \pi \in \Gamma (\bigwedge ^{2}A) \)
(i.e.of bi-degree \( (2,0) \)). In this case \( X_{\pi }=\{\pi ,\cdot \} \)
is of bi-degree \( (1,-1) \), \( F_{\pi }=\exp X_{\pi } \) is the
fibre translation by \( -d\pi  \) with respect to the fibration \( \bar{h}:{\mathcal{E}}=T^{*}\Pi A^{*}\rightarrow \Pi A^{*} \),
given in coordinates by \begin{equation}
\label{eqn:translbypi}
\begin{array}{ccl}
\tilde{x}^{i} & = & x^{i}\\
\tilde{\theta }_{a} & = & \theta _{a}\\
\tilde{p}_{i} & = & p_{i}-\frac{1}{2}\frac{\partial \pi ^{ab}}{\partial x^{i}}\theta _{a}\theta _{b}\\
\tilde{\xi }^{a} & = & \xi ^{a}-\pi ^{ab}\theta _{b}
\end{array}
\end{equation}
 where \( \pi =\frac{1}{2}\pi ^{ab}(x)\theta _{a}\theta _{b} \).
The twisting of \( \Theta  \) by \( \pi  \) is given by \( \Theta _{\pi }=F^{*}_{\pi }\Theta =(\exp X_{\pi })\Theta =\mu _{\pi }+\gamma _{\pi }+\phi _{\pi }+\psi _{\pi } \)
, where \begin{equation}
\label{eqn:twistbypi}
\begin{array}{ccl}
\mu _{\pi } & = & \mu +X_{\pi }\phi =\mu +\tilde{\pi }\phi \\
\gamma _{\pi } & = & \gamma +X_{\pi }\mu +\frac{1}{2}X^{2}_{\pi }\phi =\gamma +h_{[\pi ,\cdot ]_{\mu }}+\wedge ^{2}\tilde{\pi }\phi \\
\phi _{\pi } & = & \phi \\
\psi _{\pi } & = & \psi +X_{\pi }\gamma +\frac{1}{2}X^{2}_{\pi }\mu +\frac{1}{6}X^{3}_{\pi }\phi =\psi -d_{\gamma }\pi -\frac{1}{2}[\pi ,\pi ]_{\mu }+\wedge ^{3}\tilde{\pi }\phi 
\end{array}
\end{equation}

Again \( \Theta _{\pi } \) satisfies (\ref{eqn:structure}) if \( \Theta  \)
does. In particular, for \( \phi =0 \) this automatically produces
a new quasi-Lie bialgebroid on \( (A^{*},A) \), while for \( \psi =0 \)
we get a quasi-Lie bialgebroid on \( (A,A^{*}) \) provided \( \pi  \)
obeys the twisted Maurer-Cartan equation \begin{equation}
\label{eqn:twistedMC2}
d_{\gamma }\pi +\frac{1}{2}[\pi ,\pi ]_{\mu }=\wedge ^{3}\tilde{\pi }\phi 
\end{equation}

\begin{rem}
\label{rem:r-matrix} It may be worth mentioning that a \( \pi \in \Gamma (\wedge ^{2}A) \)
generalizes {}``\( \mathbf{r} \)-matrices'' from the Lie bialgebra
theory, while \( \omega \in \Gamma (\wedge ^{2}A^{*}) \) plays the
dual role. The twisted Maurer-Cartan equation (\ref{eqn:twistedMC2})
(resp. (\ref{eqn:twistedMC1})) is a sufficient but in general \emph{not}
a necessary condition to get a quasi-Lie bialgebroid on \( (A,A^{*}) \)
(resp. \( (A^{*},A) \)), as evidenced already by the Lie bialgebra
case. Some examples where (\ref{eqn:twistedMC2}) is in fact necessary
are considered in Section \ref{sec:examples}.
\end{rem}
From the above formulas it is not difficult to deduce the expressions
for the twisted (quasi-) differentials and brackets. For \( \Theta _{\omega } \),
\begin{equation}
\label{eqn:bracketsanddiffsomega}
\begin{array}{ccl}
d_{\mu _{\omega }} & = & d_{\mu }+[\omega ,\cdot ]_{\gamma }+\iota _{\wedge ^{2}\tilde{\omega }\psi }\\
d_{\gamma _{\omega }} & = & d_{\gamma }+\iota _{\tilde{\omega }\psi }\\
{[\alpha ,\beta ]}_{\gamma _{\omega }} & = & [\alpha ,\beta ]_{\gamma }+\tilde{\omega }\psi (\alpha ,\beta )\\
{[X,Y]}_{\mu _{\omega }} & = & [X,Y]_{\mu }+[X,Y]_{\gamma ,\omega }+\wedge ^{2}\tilde{\omega }\psi (X,Y)
\end{array}
\end{equation}
 where \( \iota  \) denotes the contraction operator, while \[
[X,Y]_{\gamma ,\omega }={\mathcal{L}}_{\tilde{\omega }X}^{\gamma }Y-{\mathcal{L}}_{\tilde{\omega }Y}^{\gamma }X-d_{\gamma }(\omega (X,Y))\]
 is \emph{}a version of the Koszul bracket. Here \( X,Y\in \Gamma (A) \),
\( \alpha ,\beta \in \Gamma (A^{*}) \) and \( {\mathcal{L}}_{\alpha }^{\gamma }=[\iota _{\alpha },d_{\gamma }]=\iota _{\alpha }d_{\gamma }+d_{\gamma }\iota _{\alpha } \).
Similarly, for \( \Theta _{\pi } \), \begin{equation}
\label{eqn:bracketsanddiffspi}
\begin{array}{ccl}
d_{\mu _{\pi }} & = & d_{\mu }+\iota _{\tilde{\pi }\phi }\\
d_{\gamma _{\pi }} & = & d_{\gamma }+[\pi ,\cdot ]_{\mu }+\iota _{\wedge ^{2}\tilde{\pi }\phi }\\
{[\alpha ,\beta ]}_{\gamma _{\pi }} & = & [\alpha ,\beta ]_{\gamma }+[\alpha ,\beta ]_{\mu ,\pi }+\wedge ^{2}\tilde{\pi }\phi (\alpha ,\beta )\\
{[X,Y]}_{\mu _{\pi }} & = & [X,Y]_{\mu }+\tilde{\pi }\phi (X,Y)
\end{array}
\end{equation}
 where \[
[\alpha ,\beta ]_{\mu ,\pi }={\mathcal{L}}_{\tilde{\pi }\alpha }^{\mu }\beta -{\mathcal{L}}_{\tilde{\pi }\beta }^{\mu }\alpha -d_{\mu }(\pi (\alpha ,\beta ))\]
 is a Koszul bracket (compare with (\ref{eqn:koszulbracket})).

\begin{rem}
The \emph{}twisting transformations \( F_{\omega } \) (resp. \( F_{\pi } \))
induce the change of splitting of \( E=A\oplus A^{*} \) by the graph
of \( \omega  \) (resp. \( \pi  \)), as can be seen from the formulas
(\ref{eqn:translbyomega}) (resp. (\ref{eqn:translbypi})). This procedure
produces new (proto-, quasi-) Lie bialgebroids whose double remains
isomorphic to the original one. The twisted Maurer-Cartan equation
(\ref{eqn:twistedMC1}) (resp. (\ref{eqn:twistedMC2})) is the condition
for the graph of \( \omega  \) (resp. \( \pi  \)) to be a Dirac
structure with respect to the Courant algebroid given by \( \Theta  \).
In the Lie bialgebroid case (\( \phi =\psi =0 \)) we recover the
strict Maurer-Cartan equation appearing in \cite{LWX1}.
\end{rem}

\section{\label{sec:symmetries}Symmetries, gauge transformations and deformations.}

The full algebra of symmetries of \( {\mathcal{E}} \) (canonical
transformations preserving the total weight) is the Lie algebra \( {\mathcal{C}^{2}}={\mathcal{C}^{0,2}}\oplus {\mathcal{C}^{1,1}}\oplus {\mathcal{C}^{2,0}} \)
of quadratic hamiltonians, acting via Poisson brackets. It is easy
to see that this action preserves the subalgebra \( {\bar{\mathcal{C}}}^{\cdot } \),
hence \( {\mathcal{C}}^{2} \) acts on the pseudo-Euclidean vector
bundle \( E=A\oplus A^{*} \). In fact, \( {\mathcal{C}}^{2} \) is
isomorphic to the \emph{Atiyah algebra} of \( E \), consisting of
infinitesimal bundle transformations preserving the canonical pairing
and covering vector fields on \( M \): 
\[
0\rightarrow {\bar{\mathcal{C}}}^{2}\rightarrow {\mathcal{C}}^{2}\rightarrow Vect(M)\rightarrow 0
\]
 where \( {\bar{\mathcal{C}}}^{2}={\mathcal{C}^{0,2}}\oplus {\bar{\mathcal{C}}}^{1,1}\oplus {\mathcal{C}^{2,0}} \)
is isomorphic to the Lie algebra of endomorphisms of \( E \) preserving
the pairing. 

The structure of \( {\bar{\mathcal{C}}}^{2} \) is quite transparent:
\( {\mathcal{C}^{0,2}}=\Gamma (\wedge ^{2}A^{*}) \) and \( {\mathcal{C}^{2,0}}=\Gamma (\wedge ^{2}A) \)
are abelian subalgebras, acted upon by \( {\bar{\mathcal{C}}}^{1,1}=\Gamma (A\otimes A^{*})=End(A) \)
in the standard way. The bracket of \( {\mathcal{C}^{2,0}} \) and
\( {\mathcal{C}^{0,2}} \) ends up in \( {\bar{\mathcal{C}}}^{1,1} \):
\( \{\pi ,\omega \}=\tilde{\pi }\tilde{\omega } \), viewed as an
operator acting on \( A \).

The group \( {\mathcal{G}} \) corresponding to \( {\bar{\mathcal{C}}}^{2} \)
will be referred to as the \emph{gauge group}; it consists of sections
over \( M \) of the bundle of Lie groups whose fibre is isomorphic
to \( SO(n,n) \), where \( n \) is the rank of \( A \). From the
above Lie algebra decomposition one expects a factorization of \( {\mathcal{G}} \)
into the product \( \Gamma (\wedge ^{2}A^{*})\times Aut(A)\times \Gamma (\wedge ^{2}A) \),
but actually this is valid only on some open subset. Of particular
interest is the factorization of the product of a \( \pi \in \Gamma (\wedge ^{2}A) \)
and an \( \omega \in \Gamma (\wedge ^{2}A^{*}) \) (both viewed as
elements of \( {\mathcal{G}} \)), which will be used in the next
section: 
\begin{equation}
\label{eqn:factorization}
\pi \omega =(\tau _{\pi }\omega )T^{-1}(\pi ,\omega )(\tau _{\omega }\pi )
\end{equation}
 This factorization is valid if and only if \( T(\pi ,\omega )=1+\tilde{\pi }\tilde{\omega } \)
is invertible, in which case 
\[
\widetilde{\tau _{\pi }\omega }=\tilde{\omega }T^{-1}=\tilde{\omega 
}(1+\tilde{\pi }\tilde{\omega })^{-1},
\]
while 
\[ 
\widetilde{\tau _{\omega }\pi }=\tilde{\pi }(T^{t})^{-1}=\tilde{\pi 
}(1+\tilde{\omega }\tilde{\pi })^{-1}. 
\]
One easily checks that \( \tau _{0}\pi =\pi  \) and \( \tau _{\omega _{1}}\tau _{\omega _{2}}\pi =\tau _{\omega _{1}+\omega _{2}}\pi  \)
(whenever the terms are defined), so one can speak of a local action
of the additive group \( \Gamma (\wedge ^{2}A^{*}) \) on the space
\( \Gamma (\wedge ^{2}A) \). This action is clearly nonlinear; in
fact, it is generated by the infinitesimal action of the abelian Lie
algebra \( {\mathcal{C}^{0,2}} \) by quadratic vector fields: \( \delta _{\omega }\pi =-\tilde{\pi }\tilde{\omega }\tilde{\pi } \).
This Lie algebra action does not integrate to a global group action,
but the corresponding transformation Lie algebroid does integrate
to a global groupoid; this groupoid induces an equivalence relation
on the space \( \Gamma (\wedge ^{2}A) \). Similarly, one has a local
action of \( \Gamma (\wedge ^{2}A) \) on \( \Gamma (\wedge ^{2}A^{*}) \).
All this is reminiscent of dressing actions of Poisson-Lie groups.

The gauge group \( {\mathcal{G}} \) acts on \( {\mathcal{E}} \)
and various objects that live there. In particular, it acts on the
space of graded Lagrangian submanifolds, preserving the subspace consisting
of those that contain \( M \) and so correspond (under the projection
\( p:{\mathcal{E}}\rightarrow \Pi E \)) to maximally isotropic subbundles
of \( A\oplus A^{*} \). For instance, given \( \omega ,\omega '\in \Gamma (\wedge ^{2}A^{*}) \),
\( \pi \in \Gamma (\wedge ^{2}A) \), we have \( F_{\omega }L_{\omega '}=L_{\omega +\omega '} \),
while \( F_{\omega }L^{*}_{\pi }=L_{\tau _{\omega }\pi }^{*} \) provided
\( \tau _{\omega }\pi  \) is defined (otherwise \( F_{\omega }L^{*}_{\pi } \)
is not \( \bar{h} \)-projectable). \( F_{\pi } \) acts in a similar
fashion.

The subalgebra preserving \( L^{*} \) (corresponding to \( A^{*}\subset A\oplus A^{*} \))
is \( {\mathcal{C}^{0,2}}\oplus {\bar{\mathcal{C}}}^{1,1} \), a semidirect
product via the standard action of \( End(A) \) on \( \Gamma (\wedge ^{2}A^{*}) \).
Therefore, the corresponding subgroup of \( {\mathcal{G}} \) acts
on Manin pairs \( (({\mathcal{E}},\Theta ),L^{*}) \) (quasi-Lie bialgebroids).
The subgroup \( \Gamma (\wedge ^{2}A^{*}) \) fixes every point of
\( L^{*} \); its action on quasi-Lie bialgebroids by twisting was
described in the previous section (formulas (\ref{eqn:twistbyomega}),
with \( \psi =0 \)). 

On the other hand, the subspace \( {\mathcal{C}^{2,0}} \) plays quite
a different role in the theory of quasi-Lie bialgebroids. Recall that
a quasi-Lie bialgebroid with \( \Theta =\mu +\gamma +\phi  \) gives
rise to a quasi-Gerstenhaber algebra structure on \( \Gamma (\wedge ^{\cdot }A) \);
in particular, one has a homotopy Lie algebra on \( \Gamma (\wedge ^{\cdot }A)[1] \)
with a differential \( l_{1}=d=d_{\gamma }=\{\gamma ,\cdot \} \)
of degree \( +1 \), a bilinear bracket \( l_{2}=[\cdot ,\cdot ]=[\cdot ,\cdot ]_{\mu }=\{\{\cdot ,\mu \},\cdot \} \)
of degree \( 0 \), and a trilinear bracket \( l_{3} \), given by
\( \phi  \) as a {}``higher derived bracket'' \( [\cdot ,\cdot ,\cdot ]=-\{\cdot ,\{\cdot ,\{\cdot ,\phi \}\}\} \),
of degree \( -1 \). One can consider the deformation theory governed
by this \( L_{\infty } \)-algebra. The appropriate structures are
elements of degree \( 1 \), \( \pi \in \Gamma (\wedge ^{2}A) \),
obeying the structure equation \begin{equation}
\label{eqn:qLAstructure}
d\pi +\frac{1}{2}[\pi ,\pi ]+\frac{1}{6}[\pi ,\pi ,\pi ]=0
\end{equation}
 It is easy to see that the last term is equal to \( -\wedge ^{3}\tilde{\pi }\phi  \),
so the above equation is exactly the same as the twisted Maurer-Cartan
equation (\ref{eqn:twistedMC2}). So each such structure twists the
quasi-Lie bialgebroid to a new one given by the formulas (\ref{eqn:twistbypi})
(with \( \psi =0 \)). The new quasi-Gerstenhaber algebra thus obtained
is given by \( d+[\pi ,\cdot]+\frac{1}{2}[\pi ,\pi ,\cdot ] \), \( [\cdot 
,\cdot ]+[\pi ,\cdot ,\cdot ] \),
\( [\cdot ,\cdot ,\cdot ] \).

\section{\label{sec:examples}Examples.}

\subsection{Arbitrary bivector fields.}

Let \( M \) be a manifold, \( A=TM \), \( {\mathcal{E}}=T^{*}\Pi TM \).
Let \( \Theta _{0}=\mu =h_{d} \), where \( d \) is the de Rham vector
field on \( \Pi TM \), and \( h_{d} \) is defined as in Section
\ref{sec:twisting}. Of course, it obeys the structure equation (\ref{eqn:structure}),
hence defines a Lie bialgebroid (with \( \gamma =0 \)). Let \( \Theta _{\pi }=F^{*}_{\pi }\Theta _{0} \)
be the twist by a bivector field \( \pi =\frac{1}{2}\pi ^{ij}\theta _{i}\theta _{j} \).
Then \( \Theta _{\pi }=\mu _{\pi }+\gamma _{\pi }+\psi _{\pi } \),
and the formulas (\ref{eqn:twistbypi}) reduce to \[
\begin{array}{ccl}
\mu _{\pi } & = & h_{d}\\
\gamma _{\pi } & = & h_{[\pi ,\cdot ]}\\
\psi _{\pi } & = & -\frac{1}{2}[\pi ,\pi ]
\end{array}\]
 Here \( [\cdot ,\cdot ]=[\cdot ,\cdot ]_{\mu } \) is the ordinary
Schouten bracket of multivector fields. Since \( \Theta _{\pi } \)
obeys (\ref{eqn:structure}), we get a quasi-Lie bialgebroid on \( (T^{*}M,TM) \).
The Maurer-Cartan equation (\ref{eqn:twistedMC2}) reduces to \( [\pi ,\pi ]=0 \),
in which case we get the Lie bialgebroid of a Poisson manifold.

Bivector fields that do not satisfy any integrability condition are
not interesting objects unless one imposes an additional structure
such as a group action satisfying various types of comapatibility
conditions (\cite{KS3},\cite{AKSM}).

\subsection{Twisted Poisson manifolds with a 3-form background.}

For \( \mu =h_{d} \) as above and a 3-form \( \phi  \), define \( \Theta _{\phi }=\mu +\phi  \).
It is immediate from (\ref{eqn:quasiLie}) that \( \Theta _{\phi } \)
satisfies (\ref{eqn:structure}) if and only if \( \phi  \) is closed.
Thus, for a closed \( \phi  \), one gets a quasi-Lie bialgebroid
structure on \( (TM,T^{*}M) \) (with \( \gamma =0 \)). Now, for
a bivector field \( \pi  \), let \( \Theta _{\phi ,\pi }=F^{*}_{\pi }\Theta _{\phi } \).
We have \( \Theta _{\phi ,\pi }=\mu _{\pi }+\gamma _{\pi }+\phi _{\pi }+\psi _{\pi } \),
where by (\ref{eqn:twistbypi}) we have:\[
\begin{array}{ccl}
\mu _{\pi } & = & h_{d}+\tilde{\pi }\phi \\
\gamma _{\pi } & = & h_{[\pi ,\cdot ]}+\wedge ^{2}\tilde{\pi }\phi \\
\phi _{\pi } & = & \phi \\
\psi _{\pi } & = & -\frac{1}{2}[\pi ,\pi ]+\wedge ^{3}\tilde{\pi }\phi 
\end{array}\]
 This defines a quasi-Lie bialgebroid structure on \( (TM,T^{*}M) \)
if and only if the twisted Maurer-Cartan equation (\ref{eqn:twistedMC2})
holds, which in this case reduces to (\ref{eqn:twistedPoisson}):\[
\frac{1}{2}[\pi ,\pi ]=\wedge ^{3}\tilde{\pi }\phi \]
This can also be thought of as the structure equation (\ref{eqn:qLAstructure})
in the quasi-Gerstenhaber algebra of polyvector fields given by the
zero differential, the Schouten bracket \( [\cdot ,\cdot ] \) and
the triple bracket \( [\cdot ,\cdot ,\cdot ]=-\{\cdot ,\{\cdot ,\{\cdot ,\phi \}\}\} \):
\[
\frac{1}{2}[\pi ,\pi ]+\frac{1}{6}[\pi ,\pi ,\pi ]=0\]
This condition can be viewed as a {}``twisted'' version of the Jacobi
identity: defining \( \{f,g\}=[[f,\pi ],g] \) and \( X_{f}=\{f,\cdot \} \)
for \( f,g\in C^{\infty }(M) \), the above equation translates to:
\[
\{\{f,g\},h\}+\{\{g,h\},f\}+\{\{h,f\},g\}=\phi (X_{f},X_{g},X_{h})\]
The triple \( (M,\pi ,\phi ) \) satisfying (\ref{eqn:twistedPoisson})
is called a \emph{\( \phi  \)-twisted Poisson manifold}. It follows
that \( (TM,T^{*}M,\mu _{\pi },\gamma _{\pi },\phi ) \) is a quasi-Lie
bialgebroid which is a deformation of the one given by \( \pi =0 \).
The deformed quasi-Gerstenhaber algebra of multivector fields is easily
seen to be given by \( d_{\gamma _{\pi }}=[\pi ,\cdot ]+\frac{1}{2}[\pi ,\pi ,\cdot ] \),
\( [\cdot ,\cdot ]_{\mu _{\pi }}=[\cdot ,\cdot ]+[\pi ,\cdot ,\cdot ] \),
with \( [\cdot ,\cdot ,\cdot ] \) unchanged. Alternatively, we can
specialize the formulas (\ref{eqn:bracketsanddiffspi}) to the present
example. The differential is given by \[
d_{\gamma _{\pi }}=[\pi ,\cdot ]+\iota _{\wedge ^{2}\tilde{\pi }\phi }\]
 where \( \iota  \) denotes contraction with the bivector-valued
1-form \( \wedge ^{2}\tilde{\pi }\phi  \). The quasi-Gerstenhaber
bracket \( [\cdot ,\cdot ]_{\mu _{\pi }} \) is uniquely determined
by\[
\begin{array}{ccl}
{[X,f]}_{\mu _{\pi }} & = & [X,f]=Xf,\\
{[X,Y]}_{\mu _{\pi }} & = & [X,Y]+\tilde{\pi }\phi (X,Y)
\end{array}\]
 where \( X \) and \( Y \) are vector fields and \( f \) is a function
on \( M \). \( d_{\gamma _{\pi }} \) squares to zero and acts as
a derivation of \( [\cdot ,\cdot ]_{\mu _{\pi }} \). The bracket
\( [\cdot ,\cdot ]_{\mu _{\pi }} \) satisfies the graded Jacobi identity
up to a homotopy given by \( \phi  \) (see (\ref{eqn:L2L2anchors})
and (\ref{eqn:L2L2brackets})); in addition, the coherence condition
(\ref{eqn:L2L3}) between \( [\cdot ,\cdot ]_{\mu _{\pi }} \) and
\( \phi  \) holds. 

This new quasi-Gerstenhaber algebra governs deformations of \( \pi  \)
within the class of \( \phi  \)-Poisson structures: \( \pi '=\pi +\delta  \)
is \( \phi  \)-Poisson if and only if \[
d_{\gamma _{\pi }}\delta +\frac{1}{2}[\delta ,\delta ]_{\mu _{\pi }}+\frac{1}{6}[\delta ,\delta ,\delta ]=0\]
 In particular, notice that \( \delta =\pi  \) does \emph{not} obey
this, since for \( \phi \neq 0 \) the equation (\ref{eqn:twistedPoisson})
is not homogeneous in \( \pi  \); furthermore, \( \pi  \) is not
a cocycle with respect to the modified differential, nor is \( [\pi ,\pi ]_{\mu _{\pi }}=0 \).

The dual quasi-differential Gerstenhaber algebra of differential forms
on \( M \) consists of the quasi-differential \[
d_{\mu _{\pi }}=d+\iota _{\tilde{\pi }\phi }\]
 where \( \iota  \) denotes the contration with the 2-form-valued
vector field \( \tilde{\pi }\phi  \), while the Schouten bracket
\( [\cdot ,\cdot ]_{\gamma _{\pi }} \) is uniquely determined (see
\ref{eqn:bracketsanddiffspi}) by \[
\begin{array}{ccl}
{[\alpha ,f]}_{\gamma _{\pi }} & = & (\tilde{\pi }\alpha )f\\
{[\alpha ,\beta ]}_{\gamma _{\pi }} & = & {\mathcal{L}}_{\tilde{\pi }\alpha }\beta -{\mathcal{L}}_{\tilde{\pi }\beta }\alpha -d_{\mu }(\pi (\alpha ,\beta ))+\wedge ^{2}\tilde{\pi }\phi (\alpha ,\beta )
\end{array}\]
 The bracket \( [\cdot ,\cdot ]_{\gamma _{\pi }} \) satisfies the
graded Jacobi identity, and \( d_{\mu _{\pi }} \) acts on it by dervations.
In addition, \( d_{\mu _{\pi }} \) squares to \( -[\phi ,\cdot ]_{\gamma _{\pi }} \),
and \( d_{\mu _{\pi }}\phi =0 \).

\subsection{Gauge transformations. }

The gauge transformations of twisted Poisson manifolds introduced
in \cite{SevWe} can be expressed in the present setting in terms
of the factorization of the gauge group described in the previous
section. Indeed, from (\ref{eqn:factorization}) we immediately get
\[
F^{*}_{\tau _{-\omega }\pi }F^{*}_{\omega }=F^{*}_{T^{-1}(-\pi ,\omega )}F^{*}_{\tau _{-\pi }\omega }F^{*}_{\pi }\]
 for a given bivector \( \pi  \) and a 2-form \( \omega  \) such
that \( T(-\pi ,\omega )=1-\tilde{\pi }\tilde{\omega } \) is invertible.
Applying both sides to \( \Theta _{\phi } \) and using formulas (\ref{eqn:twistbyomega})
and (\ref{eqn:twistbypi}), we get \begin{equation}
\label{eqn:gauge}
\Theta _{\phi -d\omega ,\tau _{-\omega }\pi }=\Phi ^{*}\Theta _{\phi ,\pi }
\end{equation}
 where \( \Phi ^{*}=F^{*}_{T^{-1}(-\pi ,\omega )}F^{*}_{\tau _{-\pi }\omega } \)
is an element of the gauge subgroup preserving \( L^{*} \). Hence,
if \( \pi  \) is a \( \phi  \)-Poisson structure, \( \tau _{-\omega }\pi  \)
is a \( (\phi -d\omega ) \)-Poisson structure. We thus have a local
action of the abelian group of 2-forms on the space of all twisted
Poisson structures, or a global transformation groupoid inducing an
equivalence relation. Equation (\ref{eqn:gauge}) implies that gauge-equivalent
twisted Poisson manifolds have isomorphic quasi-Lie bialgebroids.
Therefore, not only their Poisson cohomology spaces, but the entire
spectral sequences are isomorphic. 

The cohomology class of \( \phi  \) is preserved by gauge transformations;
locally, every \( \phi  \)-Poisson structure is equivalent to an
ordinary Poisson structure. The group of closed 2-forms acts locally
on these. If \( \pi  \) is a Poisson structure and \( \omega  \)
is a closed 2-form such that the gauge transformation \( \tau _{-\omega }\pi  \)
exists, then by (\ref{eqn:gauge}) and (\ref{eqn:twistbyomega}) \( \omega '=\tau _{-\pi }\omega  \)
obeys the Maurer-Cartan equation \( d\omega '+\frac{1}{2}[\omega ',\omega ']_{\pi }=0 \),
where \( [\cdot ,\cdot ]_{\pi } \) is the Koszul bracket. This gives
meaning to the MC equation in the differential graded Lie algebra
of differential forms on \( M \) corresponding to a Poisson structure.


\begin{thebibliography}{10}

\bibitem{AKSM}
A.~Alekseev, Y.~Kosmann-Schwarzbach, and E.~Meinrenken.
\newblock {Quasi-Poisson} manifolds.
\newblock {\em Canadian J. Math.}, 54(1):3--29, 2002.
\newblock math.DG/0006168.

\bibitem{CoSchia}
L.~Cornalba and P.~Schiappa.
\newblock Nonassociative star product deformations for {D-brane} worldvolumes
  in curved backgrounds.
\newblock preprint hep-th/0101219, 2001.

\bibitem{Dr2}
V.G. Drinfeld.
\newblock Quasi-{Hopf} algebras.
\newblock {\em Leningrad Math. J.}, 1(6):1419--1457, 1990.

\bibitem{Huebsch-Quasi}
J.~Huebschmann.
\newblock Quasi-{Lie-Rinehart} algebras, quasi-{Gerstenhaber} algebras and
  quasi-{BV} algebras.
\newblock in preparation, 2001.

\bibitem{KlimStrobl}
C.~Klim\v{c}ik and T.~Strobl.
\newblock {WZW-Poisson} manifolds.
\newblock preprint math.SG/0104189, 2001.

\bibitem{KonUrb}
K.~Konieczna and P.~Urbanski.
\newblock Double vector bundles and duality.
\newblock {\em Arch. Math. (Brno)}, 35(1):59--95, 1999.
\newblock dg-ga/9710014.

\bibitem{KS3}
Y.~Kosmann-Schwarzbach.
\newblock Jacobian quasi-bialgebras and quasi-{Poisson} {Lie} groups.
\newblock In {\em Mathematical aspects of classical field theory}, pages
  459--489. Contemp. Math., 132, Amer. Math. Soc., 1992.

\bibitem{KS5}
Y.~Kosmann-Schwarzbach.
\newblock From {P}oisson algebras to {G}erstenhaber algebras.
\newblock {\em Ann. Inst. Fourier, Grenoble}, 46(5):1243--1274, 1996.

\bibitem{SHLA}
T.~Lada and M.~Markl.
\newblock Strongly homotopy {Lie} algebras.
\newblock {\em Communications in algebra}, 23(6):2147--2161, 1995.

\bibitem{LecRog}
P.~Lecomte and C.~Roger.
\newblock Modules et cohomologie des big\'{e}bres de {Lie}.
\newblock {\em Comptes rendus {Acad. Sci. Paris}}, 310:405--410, 1990.

\bibitem{LWX1}
Zhang-Ju Liu, Alan Weinstein, and Ping Xu.
\newblock Manin triples for {Lie} bialgebroids.
\newblock {\em J. Diff. Geom.}, 45:547--574, 1997.

\bibitem{Mac-Double1}
K.~Mackenzie.
\newblock Double {Lie} algebroids and second-order geometry {I}.
\newblock {\em Adv. Math.}, 94(2):180--239, 1992.

\bibitem{MacXu}
K.C.H. Mackenzie and P.~Xu.
\newblock Lie bialgebroids and {Poisson} groupoids.
\newblock {\em Duke Math. J.}, 73:415--452, 1994.

\bibitem{Park}
J.-S. Park.
\newblock Topological open $p$-branes.
\newblock preprint hep-th/0012141, 2000.

\bibitem{Prad-1}
J.~Pradines.
\newblock G\'{e}om\'{e}trie differentielle au-dessus d'un groupo\"{i}de.
\newblock {\em C. R. Acad. Sci. Paris, s\'{e}rie A}, 266:1194--1196, 1968.

\bibitem{Roy1}
D.~Roytenberg.
\newblock {\em Courant algebroids, derived brackets and even symplectic
  supermanifolds}.
\newblock PhD thesis, UC {Berkeley}, 1999.
\newblock math.DG/9910078.

\bibitem{SevWe}
P.~\v{S}evera and A.~Weinstein.
\newblock Poisson geometry with a 3-form background.
\newblock {\em Prog. Theor. Phys. Suppl.}, 144, 2001.
\newblock math.SG/0107133.

\end{thebibliography}

\end{document}